\documentclass[a4paper, 11pt]{article}
\usepackage{amsmath, amsfonts,amssymb}
\usepackage[english]{babel}
\begin {document}

\begin{center}
\textbf {ON $C$-PROPERTIES OF THE SPACE OF IDEMPOTENT PROBABILITY MEASURES} \\

\medskip
\textbf {Ishmetov, Azad Yangibayevich$^{1}$}\\
\smallskip

{$^{1}$ Tashkent Institute of Architecture and Civil Engineering,\\ Department of Mathematics and Natural Disciplines}\\
{ishmetov\_azadbek@mail.ru}\\
\end{center}

\begin{abstract}
In the work it is shown that the space of idempotent probability measures with compact supports is kappa-metrizable if the given Tychonoff space is kappa-metrizable. It is constructed a series of max-plus-convex subfunctors of the functor of idempotent probability measures with compact supports. Further, it is established that the functor of idempotent probability measures with the compact supports preserves openness of continuous maps.
\\


\end{abstract}

The theory of idempotent measures belongs to idempotent mathematics, i. e. the field of the mathematics based on replacement of usual arithmetic operations with idempotent (as, for example, $x\oplus y=\mbox {max}\left\{x,y\right\}$). The idempotent mathematics intensively develops at this time (see, for example, [1], survey article [2] and the bibliography in it). Its communication with traditional mathematics is described by the informal principle according to which there is a heuristic compliance between important, interesting and useful designs the last and similar results of idempotent mathematics.

In the present article we investigate a functor $I_{\beta}$  which is an extension of the functor of idempotent probability measures from the category of compact Hausdorff spaces onto the category of Tychonoff spaces and their continuous maps. In traditional mathematics to it there corresponds the functor $P_{\beta}$  of probability measures. The concept of an idempotent measure (Maslov's measure) finds numerous applications in various field of mathematics, mathematical physics and economy. In particular, such measures arise in problems of dynamic optimization; the analogy between Maslov's integration and optimization is noted also in [1]. It is well-known that use of measures of Maslov for modeling of uncertainty in mathematical economy can be so relevant as far as also use of classical probability theory.

Unlike a case of probability measures to which consideration extensive literature is devoted topological properties of the spaces of idempotent measures were practically not investigated. In work [2] M.Zarichnyi gave a number of appendices of idempotent measures in various branches of modern sciences.

Let $\mathbb{R}=\left( -\infty,+ \infty \right)$  be the real line. On the set $\mathbb{R}\bigcup \left\{ -\infty  \right\}$  we define operations  $\oplus $  and $\odot$ by the rules: $u\oplus v=\max \{u,v\}$  and  $u\odot v=u+v$. It is easy to see $-\infty$  is the zero  $\mathbf{0}$, and the usual zero $0$  is the unit $\mathbf{1}$  on  $\mathbb{R}\bigcup \left\{ -\infty  \right\}$. The collection $\left(\mathbb{R}\bigcup \left\{ -\infty \right\}, \oplus, \odot,\mathbf{0}, \mathbf{1} \right)$,  forms the `max-plus' semi-field which we denote by  $\mathbb{R}_{\max}$.

Let $X$  be a compact (i. e. Hausdorff compact space, $\it {plural:}$ compacts. Note that a compactum is a metrizable compact space, $\it {plural:}$ compacta), $C(X)$  be the algebra of all continuous functions defined on  $X$.  $C(X)$ is endowed with the usual pointwise algebraic operations and  $\sup$-norm. We introduce the following operations:

1)  $\odot:\mathbb{R}\times C(X)\rightarrow C(X)$ by a rule $\odot
(\lambda, \varphi)=\lambda\odot\varphi=\varphi+\lambda_{X}$, where $\varphi\in C(X)$  and  $\lambda_{X}$ is constant function accepting everywhere on $X$  the value  $\lambda\in \mathbb{R}$;

2) $\oplus: C(X)\times C(X) \rightarrow C(X)$  by a rule èëó  $\oplus
(\varphi, \psi)=\varphi \oplus\psi = \max\{\varphi,\psi\}$, where $\varphi,\psi\in C(X)$.

{\bf Definition 1[2].} A functional $\mu :C(X)\to \mathbb{R}$ is called {\it an idempotent probability measure on} $X$ if
it satisfies the following properties:

(i) $\mu ({{\lambda }_{X}})=\lambda $ for any $\lambda \in \mathbb{R}$ ({\it norm axiom});

(ii) $\mu(\lambda \odot \varphi)=\lambda \odot \mu(\varphi)$ for any $\lambda \in \mathbb{R}$ and $\varphi \in C(X)$ ({\it homogeneity axiom});

(iii) $\mu (\varphi \oplus \psi )=\mu (\varphi )\oplus \mu (\psi )$ for any $\varphi $, $\psi \in C(X)$ ({\it additivity axiom}).

The number $\mu(\varphi)$  is called the Maslov's integral corresponding to  $\mu$. The set of all idempotent probability measure on $X$  we denoted by  $I(X)$. We have  $I(X)\subset \mathbb{R}^{C(X)}$. Consider $I(X)$  with induced from  $\mathbb{R}^{C(X)}$ topology. The sets of the look
$$
\left\langle\mu;\ \varphi_1,\ ...,\ \varphi_k;\ \varepsilon\right\rangle = \{\nu\in I(X):\ |\nu(\varphi_i)-\mu(\varphi_i)|<\varepsilon,\ i=1,\ ...,\ k\}
$$

form a base of neighborhoods of an idempotent probability measure $\mu \in I(X)$  concerning to this topology. form a base of neighborhoods of an idempotent probability measure   concerning to this topology. Here  $\varphi_1,...,\varphi_k\in C(X)$ and $\varepsilon>0$. It is well known that for any compact $X$  the space $I(X)$  is also a compact. Let $f:X\rightarrow Y$  be a continuous map of compacts. Then the equality
 $$I(f)(\mu)(\varphi )=\mu (\varphi \circ f), \ \ \mu\in I(X),\ \  \varphi\in C(Y),$$
defines a map $I(f):I(X)\to I(Y)$  which is continuous. For an idempotent probability measure  $\mu \in I(X)$ we define  its support:
 $$
\text{supp}\mu =\bigcap \left\{ F\subset X: F \mbox{ is a closed subset of } X  \mbox{ and } \mu \in I(F) \right\}.
$$
For a compact $X$  and  positive integer $n$ we put
 $${{I}_{n}}(X)=\left\{ \mu \in I(X):\,\,\left| \operatorname{supp}\mu  \right|\le n \right\}.$$
Further
 $${{I}_{\omega }}(X)=\bigcup\limits_{n=1}^{\infty }{{{I}_{n}}(X)}.$$

Let $X$  is Tychonoff space, $\beta X$  be the Stone-\v{C}ech compact extension of $X$. We define [11, 15] a subspace
 $$I_{\beta}(X)=\{\mu\in I(\beta X): \mbox{supp}\mu\subset X\},$$
which elements we call as {\it idempotent probability measures with compact support}. Let $\beta f: \beta X \rightarrow \beta Y$, where be the maximal extension of a continuous map $f:X\rightarrow Y$  of Tychonoff spaces. Then $I(\beta f)(I_{\beta}(X))\subset I_{\beta}(Y)$. Put $$I_{\beta}(f)=I(\beta f)\mid I_{\beta}(X).$$

Thus, the operation $I_{\beta}$  is a functor acting in the category $Tych$  Tychonoff spaces and their continuous maps.

For positive integer $n$  put  $I_{\beta, n}(X)=\left\{\mu\in I_{\beta}(X): \,\,\left|\operatorname{supp}\mu  \right|\le n\right\}$. Put  ${{I}_{\beta, \omega }}(X)=\bigcup\limits_{n=1}^{\infty }{{{I}_{\beta, n}}(X)}$.

{\bf Proposition 1.} If  $Y$  is everywhere dense in a compact  $X$, then $I_{\beta,\omega}(Y)$  is everywhere dense in $I(X)$.

{\bf Proof.} It is well-known [11] that $I_{\omega}(X)$  is everywhere dense in  $I(X)$. Therefore it is enough to establish that $I_{\beta,\omega}(Y)$ is everywhere dense in  $I_{\omega}(X)$. Take a measure $\mu\in I_{\omega}(X)$ and its basic neighbourhood  $\langle\mu;\varphi_1, ... , \varphi_k; \varepsilon\rangle$. Let $\mu=\lambda_1\odot\delta_{x_1}\oplus\lambda_2\odot\delta_{x_2}\oplus...\oplus\lambda_s\odot\delta_{x_s}$. As $Y$  is everywhere dense in  $X$, there are points $y_1,...,y_s$  such that $|\varphi_i(x_j)-\varphi_i(y_j)|<\varepsilon$  for all  $i=1,...,k$;  $j=1,...,s$. There exist  $\lambda_1',...,\lambda_s'$, that $|\lambda_i-\lambda_j'|<\varepsilon$  for all  $j=1,...,s$. That is why      $\nu=\lambda_1'\odot\delta_{y_1}\oplus\lambda_2'\odot\delta_{y_2}\oplus...\oplus\lambda_s'\odot\delta_{y_s}\in \langle\mu;\varphi_1, ... , \varphi_k; \varepsilon\rangle\cap I_{\beta, \omega}$. Proposition 1 is proved.

{\bf Corollary 1.} If $Y$  is an everywhere dense subspace of a compact  $X$, then $I_{\omega}(Y)$  and $I_{\beta}(Y)$ are everywhere dense subspaces of  $I(X)$.

{\bf Proposition 2.} If $Y$ is everywhere dense subspace of a Tychonoff space $X$, then $I_{\beta,\omega}(Y)$  is an everywhere dense subset of  $I_\beta(X)$.

{\bf Proof.} Let $bX$  be a compact extension of  $X$. Then  $Y$, being everywhere dense in  $X$, is everywhere dense in  $bX$. According to the proposition 1 the set  $I_{\beta,\omega}(Y)$ is everywhere dense in  $I(bX)$. But  $I_{\beta,\omega}(X)\subset I_{\beta}(X)\subset I(bX)$. Therefore,  $I_{\beta,\omega}(Y)$  is everywhere dense in  $I_{\beta}(X)$ . Proposition 2 is proved.

{\bf Corollary 2.} For every Tychonoff space $X$  the set $I_{\beta,\omega}(X)$  is everywhere dense in  $I_{\beta}(X)$.

{\bf Definition 2[3].} Let  $\mathcal{P}$ be some topological property. A Tychonoff space  $X$ is called  $C-\mathcal{P}$-space if it has a compact extension  $bX$, satisfying the property  $\mathcal{P}$.

Objects of our attention are  $C$-dyadic spaces,  $C$-Milyutin spaces,  $C$-Dugundji spaces,  $C$-absolute retracts or  $C$-$AR$-spaces (see [3]). To research of the specified classes of spaces we need the following auxiliary statement: if $X$  is a Tychonoff space of weight $\leq\tau$  and  $bX$ is its compact extension which is a dyadic compact, then $wbX\leq\tau$ [3].

Note that for a topological space $X$  its weight (i. e. the smallest power of bases of $X$) is denoted by $wX$.

{\bf Theorem 1.} If $X$  is a  $C$-dyadic space of the weight $\leq\omega_1$   then $I_{\beta}$  is also a  $C$-dyadic space.

{\bf Proof.} Let $bX$  be a compact extension of $X$  which is dyadic. Then the weight of $bX$  is not more than  $\omega_1$. Then there is an epimorphism  $f:D^{\omega_1}\rightarrow bX$. But $D^{\omega_1}$  is a Dugundji compact. Therefore,  $I_{\beta}(D^{\omega_1})$ is an absolute retract according to [18]. But every $AR$-compact is dyadic. Therefore, the compact $I_{\beta}(bX)$  is also dyadic, being image of a dyadic compact  $I_{\beta}(D^{\omega_1})$ rather continuous map. At last, by a Corollary 1 the space $I(bX)$  is a compact extension of the space  $I_{\beta}(bX)$. Hence $I_{\beta}(X)$  is  $C$-dyadic. Theorem 1 is proved.

{\bf Proposition 3[3].} For Tychonoff space $X$  of the weight $\leq\omega_1$  the following conditions are equivalent:

1)  $X$ is  $C$-Milyutin space;

2)  $X$ is  $C$-Dugundji space.

{\bf Theorem 2.} Let $X$  be a  $C$-Milyutin space of weight  $\leq\omega_1$. Then $I_{\beta}(X)$  is  $C$-absolute retract.

{\bf Proof.} According to Proposition 3 $X$  is a  $C$-Dugundji space. Take a compact extension  $bX$ of $X$, which is a Dugundji compact. Therefore  $wbX\leq\omega_1$. Then the compact $I(bX)$   is an absolute retract according to [18]. But by Corollary 1 $I(bX)$  is a compact extension of  $I_{\beta}(X)$, from here follows that $I_{\beta}(X)$  is a  $C$-$AR$-space. Theorem 2 is proved.

Since every  $AR$-compact is a Dugundji space, Theorem 2 implies

{\bf Corollary 3.} Functor $I_{\beta}$  translates the class of  $C$-absolute retracts of weight $\leq\omega_1$  into the class of  $C$-absolute retracts.

{\bf Definition 3[3].}  $\kappa$-metric (kappa-metric) on a Tychonoff space $X$  is a non-negative function $\rho(x,C)$  of two variables: points $x\in X$  and canonically closed sets  $C=[<C>]\subset X$, satisfying to the following axioms:

Ê1) ({\it belongings axiom}).	$\rho(x,C)=0$  if and only if $x\in C$;

Ê2) ({\it monotonicity axiom}).If $C' \subset C$  then  $\rho(x,C)\leq
\rho(x,C')$;

Ê3) ({\it continuity axiom}).	At fixed $C$  the function  $\rho(x,C)$ is continuous by  $x$;

Ê4) ({\it union axiom}).	$\rho(x,[\bigcup\limits_{\alpha}C_{\alpha}])=\inf\limits_{\alpha}\rho(x,C_{\alpha})$  for any increasing well-ordered sequence of canonically closed sets  $C_{\alpha}\subset X$.

{\bf Theorem 3.} If a Tychonoff space $X$  is $C$-$\kappa$-metrizable, then $I_{\beta}(X)$  is also  $C$-$\kappa$-metrizable.

{\bf Proof.} Let $bX$  be a $\kappa$-metrizable compact extension of  $X$. By E. V. Shchepin's theorem a class of the  $\kappa$-metrizable compacts coincides with a class of the open generated compacts [3]. Here, a compact is open generated if it is homeomorphic to the limit space of some countable-directed continuous inverse spectrum $S$  consisting of compacta and open projections. Let  $bX=\lim\limits_{\leftarrow}S$ be the above stated representation of the open generated compact  $bX$. Then owing to the continuity [3] of the functor $I$  we have
 . 						$$I(bX)=\lim\limits_{\leftarrow}I(S). \eqno(1)$$
On the same reason the inverse spectrum $I(S)$  is continuous. Projections of the spectrum  $I(S)$ are open [7, 8]. Thus, equality (1) gives that the compact $I(bX)$  is open generated, i. e.  $I(bX)$ is kappa-metrizable, and as it was noted above (see Corollary 1), it is a compact extension of  $I_{\beta}(X)$. Thus, $I_{\beta}(X)$  is  $C$- kappa-metrizable. Theorem 3 is proved.

While the proof of the main result we should use the following two lemmas proved in [4 -- 6].

{\bf Lemma 1.} Let $f:X\rightarrow Y$   be a continuous map, $y_0 \in Y$  and $\varphi\in C_b(Y)$.  Then there is a function $\psi\in C_b(Y)$  such that $\psi\circ f\leq \varphi$  and  $\psi(y_0)=inf\left\{\varphi(x): x\in f^{-1}(y_0)\right\}$.

{\bf  Lemma 2.} Let $f:X\rightarrow Y$  be a continuous map,  $y_0 \in Y$ and $\nu\in I_{\beta}$  such that $I_{\beta}(f)(\nu)=\delta_{y_0}$. Then for any $\varphi\in C_b(X)$  such that  $\varphi(x)\geq c$ ($\varphi(x)\leq c$) at all $x\in f^{-1}(y_0)$   we have $\nu(\varphi)\geq c$  (respectively, $\nu(\varphi)\leq c$).

Remind that a subset $A$  of $I(X)$  is max-plus-convex if $\alpha\odot\mu\oplus\beta\odot\nu\in I(X)$  for every pair of measures  $\mu, \nu\in I(X)$  where $\alpha, \beta\in \mathbb{R}_{\max}$  and  $\alpha\oplus\beta=\mathbf{1}(=0)$.

{\bf Proposition 3.} For a map $f:X\rightarrow Y$  of compacts and every measure $\nu\in I(X)$  the preimage $I(f)^{-1}(\nu)$  is a max-plus-convex set in  $I(X)$.

{\bf Proof.} Let  $\mu_1, \mu_2\in I(f)^{-1}(\nu)$. Then for all $\alpha, \beta\in \mathbb{R}_{\max}$  with $\alpha\oplus\beta=\mathbf{1}$  we have

 $$I(f)(\alpha\odot\mu_1\oplus\beta\odot\mu_2)(\psi)=(\alpha\odot\mu_1\oplus\beta\odot\mu_2)(\psi\circ f)=\alpha\odot\mu_1(\psi\circ f)\oplus\beta\odot\mu_2(\psi\circ f)=$$
 $$=\alpha\odot I(f)(\mu_1)(\psi)\oplus\beta\odot I(f)(\mu_2)(\psi)=\alpha\odot\nu(\psi)\oplus\beta\odot\nu(\psi)=\nu(\psi),$$
 $\psi\in C(Y)$.  So, $I(f)(\alpha\odot\mu_1\oplus\beta\odot\mu_2)=\nu$. Proposition 3 is proved.

{\bf Proposition 4 [16].} The set $I_{\beta}$  is a $\mbox {max}$-$\mbox{plus}$-convex subset of  $\mathbb{R}^{C_b(X)}$.

{\bf Definition 4.} A subfunctor $F$ of the functor $I_{\beta}$, acting in the category  $Tych$, is called  $\mbox {max}$-$\mbox{plus}$-convex if for any Tychonoff space $X$  the space  $F(X)$ is a  $\mbox{max}$-$\mbox{plus}$-convex subset of  $I_{\beta}(X)$.

Equivalent definition looks as follows. A subfunctor $F$  of $I_{\beta}$  is $\mbox{max}$-$\mbox{plus}$  if for any Tychonoff space  $X$, for each pair $\mu_{1},\mu_{2}\in F(X)$  and for all $\alpha, \beta\in \mathbb{R}_{\max}$, $\alpha\oplus\beta=\mathbf{1}$  we have $\alpha \odot \mu_{1}\oplus \beta \odot\mu_{2}\in F(X)$.

For a cardinal number $\tau$  we denote by $I_{\tau}$  the operation which puts in compliance to every Tychonoff space $X$  the set $I_{\tau}(X)$  of all measures $\mu\in I_{\beta}(X)$  which support's power less then  $\tau$, and to each continuous map $f:X\rightarrow Y$  a map $I_{\tau}(f)$  which is the restriction of $I_{\beta}(f)$  on $I_{\tau}(X)$.

The following statement gives a large class of  $\mbox{max}$-$\mbox{plus}$-convex subfunctors of  $I_{\beta}$.

{\bf Theorem 4.} Let $\tau$  be an infinite cardinal number. Then $I_{\tau}$  is a  $\mbox{max}$-$\mbox{plus}$-convex subfunctor of  $I_{\beta}$.

{\bf Proof.} For a cardinal number $\tau$  and a Tychonoff space  $X$ by definition we have
 $$I_{\tau}(X)=\{\mu\in I_{\beta}(X):\ \ |\mbox{supp}\mu|<\tau\}.$$
For every continuous map  $f:X\rightarrow Y$ we have  $I_{\tau}(f)=I_{\beta}(f)|_{I_{\tau}(X)}$. Since
 $$\mbox{supp}I_{\tau}(f)(\mu)=\mbox{supp} I_{\beta}(f)(\mu)=f(\mbox{supp}\mu),\ \ \ \mu\in I_{\tau}(X),$$
and
$$|\mbox{supp} I_{\tau}(f)(\mu)|=|f(\mbox{supp}\mu)|\leq |\mbox{supp}\mu|<\tau$$
we have  $I_{\tau}(f)(I_{\tau}(X))\subset I_{\tau}(Y)$, i. e.  $I_{\tau}(f)$ is a map from $I_{\tau}(X)$  into  $I_{\tau}(Y)$. The preservation of compositions of maps and identical map by the operation  $I_{\tau}$ is obvious. Therefore, $I_{\tau}$  is a subfunctor of  $I_{\beta}$.

Let's check its convexity. Let $X$ be a Tychonoff space, $\mu_{1},\mu_{2}\in I_{\beta}(X)$, $\alpha, \beta\in \mathbb{R}_{\max}$, $\alpha\oplus\beta=\mathbf{1}$  and  $\mu=\alpha\odot\mu_1\oplus\beta\odot\mu_2$. If $\beta$  or $\alpha$  equals to $-\infty$  then $\mu$  coincides with a measure $\mu_1$  or  $\mu_2$, respectively. If $\alpha >-\infty$  (in this case $\beta=0$  ), or  $\beta >-\infty$ (in this case $\alpha=0$ ) then
 $$\mbox{supp}\mu=\mbox{supp}\mu_{1}\cup\mbox{supp}\mu_{2}.$$
Anyway
 $$\mbox{supp}\mu\subset \mbox{supp}\mu_{1}\cup \mbox{supp}\mu_{2}.$$
Therefore,
 $$|\mbox{supp}\mu|\leq |(\mbox{supp}\mu)_{1}|+|\mbox{supp}\mu_{2}|.$$
From here taking into account the infinity of the cardinal number $\tau$  we receive  $|\mbox{supp}\mu|<\tau$, i. e.  $\mu\in I_{\tau}(X)$. Theorem 4 is proved.

Let $f:X\rightarrow Y$  be a continuous map and  $\varphi\in C_b(X)$. By $\varphi^{*}$  (respectively,  $\varphi_{*}$) we denote a function $\varphi^{*}:Y\rightarrow \mathbb{R}$  (respectively,  $\varphi_{*}:Y\rightarrow \mathbb{R}$) defined by the rule $\varphi^{*}(y)=\mbox{sup}\left\{\varphi(x): x\in f^{-1}(y)\right\}$  (respectively, $\varphi_{*}(y)=\mbox{inf}\left\{\varphi(x): x\in f^{-1}(y)\right\}$ ). It is known that if $f$  is an open map then the functions $\varphi^{*}$  and $\varphi_{*}$  are continuous.

{\bf Theorem 5.} Let $f:X\rightarrow Y$  be a continuous map from a Tychonoff space  $X$ to a Tychonoff space $Y$. Then the map $I_{\beta}(f):I_{\beta}(X)\rightarrow I_{\beta}(Y)$  is open if and only if  $f$ is open.

{\bf Proof.} Let $f:X\rightarrow Y$  be such a map that the map $I_{\beta}(X)\rightarrow I_{\beta}(Y)$  is open. Fix a point  $x_0\in X$. Let  $y_0=f(x_0)$. Take such  $\varphi\in C_b(X)$, that  $\varphi(x_0)=0$. Put
$$V=\left\{x\in X: -1<\varphi(x)<1\right\}.$$
As the sets of the view $V$  form a base of neighbourhood of the point $x_0$  it is sufficient to show that $f(V)$  is an open neighbourhood of  $y_0$. Consider an open neighbourhood $W=\left\{\mu\in I_{\beta}(X): -1<\mu(\varphi)<1\right\}$  of the idempotent probability measure  $\delta_{x_0}\in I_{\beta}(X)$. Then $I_{\beta}(f)(W)$  is an open neighbourhood of the idempotent probability measure  $\delta_{y_0}\in I_{\beta}(Y)$. There are functions $\psi_1, \psi_2,...,\psi_k\in C_b(Y)$  with $\psi_i(y_0)$  and  $\varepsilon$ with $0<\varepsilon<1$, such that  $H=\langle\delta_{y_0},\psi_1, \psi_2,...,\psi_k;\varepsilon\rangle\subset I_{\beta}(f)(W)$.

Put $G=\bigcap\limits^{n}_{i=1}\left\{y\in Y: -\varepsilon<\psi_{i}(y)<\varepsilon\right\}$. Then $G$  is an open neighbourhood of the point  $y_0$. Let $y\in G$  be an arbitrary point. Then  $\delta_y\in H$. Consequently, there exists $\mu\in I_{\beta}(X)$  such that $\mu\in W$  and  $I_{\beta}(f)(\mu)=\delta_y$. By the condition  $-1<\mu(\varphi)<1$. Since every idempotent probability measure is an order-preserving functional, then $\mu$  is an order-preserving functional, that is why there exists $x\in f^{-1}(y)$  such that  $-1<\varphi(x)<1$. Thus, $G\subset f(V)$   and therefore $f$  is open.

Let now $f:X\rightarrow Y$  be an open map. Assume $I_{\beta}(f)$  is not open. Then there exist:

1) an idempotent probability measure  $\mu_0\in I_{\beta}(X)$,

2) a net of idempotent probability measures $\left\{\nu_{\alpha}\right\}\subset I_{\beta}(X)$  converging to $\nu_0=I_{\beta}(f)(\mu_0)$   and

3) a neighbourhood $W$  of  $\mu_0$ such that $I_{\beta}(f)^{-1}(\nu_{\alpha})\bigcap W=\oslash$  for every $\alpha$.

As  $I_{\omega}(Y)$ is everywhere dense in $I_{\beta}(Y)$  one can assume that all of $\nu_{\alpha}$  are idempotent probability measures with finite support. Put $A_{\alpha}=I_{\beta}(X)^{-1}(\nu_{\alpha})$  and  $\left[A_{\alpha}\right]=\left[I_{\beta}(f)^{-1}(\nu_{\alpha})\right]_{I(\beta X)}$. As $I(\beta X)$  is a compact and the function $I(\beta f)$  is continuous the net $\left[A_{\alpha}\right]$  converges to $\left[A_{0}\right]=\left[I_{\beta}(f)^{-1}(\nu_{0})\right]_{I(\beta X)}$  according to Vietoris topology in  $\mbox {exp}I(\beta X)$. Besides, owing to continuity of the map $I(\beta f)$  we have $\left[A_{0}\right]\subset I(\beta f)^{-1}(\nu_{0})$  and  $\mu_0\notin\left[A_{0}\right]$ (otherwise condition 3) is violated). As $\mu_0\notin\left[A_{0}\right]$  for every $\mu\in\left[A_{0}\right]$  there exists $\varphi_{\mu}\in C_b(X)\cong C(\beta X)$  such that  $\mu_0\left(\varphi_{\mu}\right)\neq\mu\left(\varphi_{\mu}\right)$. According to the assumption for every $\alpha$  there exists a finite set $\left\{y_{\alpha 1},...,y_{\alpha {n_{\alpha}}}\right\}$  such that  $\nu_{\alpha}\in I\left(\left\{y_{\alpha 1},...,y_{\alpha {n_{\alpha}}}\right\}\right)$. For every $y_{\alpha i}$  choose $x_{\alpha i}\in X$  such that $f\left(x_{\alpha i}\right)=y_{\alpha i}$  and  $\varphi\left(x_{\alpha i}\right)=\varphi^{*}\left(y_{\alpha i}\right)$,  $\mu\in\left[A\right]$. Define an embedding $j_{\alpha}:\left\{y_{\alpha 1},...,y_{\alpha {n_{\alpha}}}\right\}\rightarrow X$  by the rule $j_{\alpha}(y_{\alpha i})=x_{\alpha i}$  and put $\mu_{\alpha}=I(j_{\alpha})(\nu_{\alpha})$. It is easy to see that $\varphi=\varphi^{*}\circ f=\varphi^{*}\circ \beta f$  on every  $\left\{y_{\alpha 1},...,y_{\alpha {n_{\alpha}}}\right\}$. That is why $$\mu_{\alpha}(\varphi)=\mu_{\alpha}\left(\varphi^{*}\circ f\right)=\mu_{\alpha}\left(\varphi^{*}\circ \beta f\right)=I(\beta f)(\mu_{\alpha})(\varphi^{*})=\nu_{\alpha}(\varphi_{*})$$

for every $\alpha$. Let $\mu_1$   be a limit of the net  $(\mu_{\alpha})$. Then $\mu_1\in\left[A\right]$  and $$\mu_{1}(\varphi)=\lim\limits_{\alpha}\mu_{\alpha}(\varphi)=\lim\limits_{\alpha}\mu_{\alpha}\left(\varphi^{*}\circ \beta f\right)=\lim\limits_{\alpha}I(\beta f)(\mu_{\alpha})(\varphi^{*})=\lim\limits_{\alpha}\nu_{\alpha}(\varphi^{*})=\nu_{0}(\varphi^{*}),\ \ \varphi\in C(X).$$

On the other hand  $\nu_0(\varphi^{*})=I(\beta f)(\mu_0)(\varphi^{*})=\mu_0(\varphi^{*}\circ\beta f)$. Thus,
 						$$\mu_1(\varphi)=\mu_0(\varphi^{*}\circ\beta f) \eqno (2)$$
for every $\varphi\in C(X)$. Similarly,
 						$$\mu_1(\varphi)=\mu_0(\varphi_{*}\circ\beta f) \eqno (3)$$
for every  $\varphi\in C(X)$. Let $\varphi_{\mu_1}\in C(X)$  be a function such that  $\mu_0(\varphi_{\mu_1})\neq\mu_1(\varphi_{\mu_1})$. Suppose  $\mu_0(\varphi_{\mu_1})>\mu_1(\varphi_{\mu_1})$. Since $\varphi^{*}\beta f\geq\varphi$  owing to (2) we have  $\mu_1(\varphi_{\mu_1})=\mu_0(\varphi^{*}_{\mu_1}\circ\beta f)\geq\mu_0(\varphi_{\mu_1})$. Analogously one can show the assumption $\mu_0(\varphi_{\mu_1})<\mu_1(\varphi_{\mu_1})$  also is false. Thus, we get a contradiction, which shows that $I_{\beta}(f)$  is open. Theorem 5 is proved.

{\bf Acknowledgement.} The author would like to thank to professor Adilbek Zaitov – the head of the department of Mathematics and Natural Disciplines of Tashkent institute of architecture and civil engineering for comprehensive support and attention.

\begin{center}
\textsl{References}
\end{center}

[1].  Maslov V. P. Operators methods. -- Moscow, `Mir'. 1987 (Russian).

[2]. Zarichnyi M. Idempotent probability measures, I. arXiv:math. GN/0608754 v 130 Aug 2006.

[3]. Shchepin E. V. Functors and uncountable powers of compacta. //Russian Math. Surveys, 1981, Vol 36 Issue 3, P. 1-71.

[4].	Zaitov A. A. On categorical properties of order-preserving functionals. //Methods of Functional Analysis and Topology. vol. 9(2003), ¹4, p. 357-364.

[5].	Zaitov A. A. On extension of order-preserving functionals. //Doklady Akademii Nauk of Uzbekistan. 2005. No 5. P. 3-7.

[6].	Zaitov A. A. On open mapping theorem of the spaces of order-preserving functionals. //Mathematical notes. Vol. 88. No 5. 2010. P. 21-26

[7].	Zaitov A. A., Tojiev I. I. Functional representations of closed subsets of a compact. //Uzbek mathematical journal, 2010. No 1. P. 53-63.

[8].	Zaitov A. A., Tojiev I. I. On uniform metrizability of the functor of idempotent probability measures. //Uzbek mathematical journal, 2011. No 2. P. 66-74. DOI: On uniform metrizability of the functor of idempotent probability measures. //arxiv: 1204. 0074v1 [math. GN] 31 March 2012.

[9].	Zaitov A. A., Kholturayev Kh. F. On perfect metrizability of the functor of idempotent probability measures. //arxiv: 1205. 0864v1 [math. GN] 4 May 2012.

[10].	Zaitov A. A., Tojiev I. I. On metric of the space of idempotent probability measures //arxiv:1006. 3902 v 2 [math. GN] 15 March 2012.

[11].	Zaitov A. A., Ishmetov A. Ya. On monad generating by the functor  . //Vestnik of National University of Uzbekistan. 2013. No 2. P. 61-64.

[12].	Zaitov A. A., Kholturaev Kh. F. On interrelation of the functors   of probability measures and   of idempotent probability measures. //Uzbek Mathematical Journal, 2014. No 4. P. 36-45.

[13].	Zaitov A. A., Ishmetov A. Ya. Geometrical properties of the space   of idempotent probability measures. //arxiv:1808. 10749v2 [math. GN] 4 Sep 2018.

[14].	Zaitov A. A., Kholturaev Kh. F. Geometrical properties of the space of idempotent probability measures. //arXiv:1811. 08325v1 [math. GN] 19 Nov 2018.

[15].	Ishmetov A. Ya. On functor of idempotent probability measures with compact support. //Uzbek mathematical journal, 2010, No 1. P. 72-80.

[16].	Ishmetov A. Ya. On max-plus-convex subfunctors of the functor   of idempotent probability measures. //Sciences of Europe, Vol 3, No 33, P. 47-54.

[17].	Tojiev I. I. On a metric of the space of idempotent probability measures. //Uzbek mathematical journal, 2010. No 4. P. 165-172.

[18].	Radul T. Idempotent measures: absolute retracts and soft maps. //arxiv: 1810.09140 v1. [math.GN] 22 Oct 2018.

\end{document}